\input amstex
\documentstyle{amsppt}

\hsize=4.75in
\vsize=8in

	
	\let	\< = \langle
	\let	\> = \rangle

	\define		\a		{\alpha}
	\redefine	\b		{\beta}
	\redefine	\d		{\delta}
	\redefine	\D		{\Delta}
	\define		\e		{\varepsilon}
	\define		\g		{\gamma}
	\define		\G		{\Gamma}
	\redefine	\l		{\lambda}
	\redefine	\L		{\Lambda}
	\define		\n		{\nabla}
	\redefine	\var	{\varphi}
	\define		\s		{\sigma}
	\redefine	\Sig	{\Sigma}
	\redefine	\t		{\tau}
	\define		\th		{\theta}
	\redefine	\O		{\Omega}
	\redefine	\o		{\omega}
	\define		\z		{\zeta}
	
	\redefine	\i		{\infty}
	\define		\p		{\partial}

\rightheadtext{The space-time uncertainty relations}
\leftheadtext{Marc A. Rieffel}
\topmatter
\title
On the operator algebra for the space-time uncertainty relations
\endtitle
\author
Marc A. Rieffel
\endauthor
\thanks The research reported here 
   was supported in part by National Science Foundation grant
DMS--96-13833.
\endthanks
\affil
Department of Mathematics \\
University of California \\
Berkeley, CA\ \ 94720
\endaffil
\email rieffel\@math.berkeley.edu
\endemail
\endtopmatter
\document

\heading
1.  Introduction
\endheading
The purpose of this note is to show that the construction of the $C^*$-algebra
for the space-time uncertainty relations which was introduced by Doplicher,
Fredenhagen and Roberts \cite{2,3,4} fits comfortably into the deformation
quantization framework developed in \cite{5}.  This has the mild advantages 
that one can work directly with functions on space-time 
rather than with their Fourier
transforms, the treatment of the unbounded space-time operators is fairly
smooth, and the sense in which one has a deformation of commutative space-time
is made technically precise.  Possibly our techniques will be useful in some
related situations.

We will not repeat here the physical motivation and the treatment of the
uncertainty relations themselves given in \cite{2,3,4}.  We will deal only with
the construction of the $C^*$-algebra and the affiliated operators having the
desired properties.  To the extent that there is no additional complication, we
will work in slightly greater generality than immediately needed, both because
this might be useful at some later time, and because it clarifies which aspects
of structure are involved.

We now recall from \cite{4} the desired mathematical properties, in the form
most convenient for our purposes.  Let $M_0$ denote Minkowski space, and let
$L$ denote the full Lorentz group.  Let $L$ act on $4 \times 4$ matrices by
similarity, that is, $\Lambda \in L$ sends the matrix $\s$ to
$\Lambda\s\Lambda^t$ (where $t =$ transpose).  Let $\s_0$ denote the standard
symplectic matrix (as in $3.28$ of \cite{4}), and let $\Sigma$ denote its orbit
under $L$.  Then $\Sigma$ is a manifold consisting of certain invertible
skew-symmetric matrices, and is a homogeneous space for $L$.  (See $3.24$ of
\cite{4}.)  We want unbounded self-adjoint operators ${\bold q}_{\mu}$, $\mu =
0,\dots,3$, corresponding to the standard coordinates on space-time.  Let
$$
Q_{\mu\nu} = -i[{\bold q}_{\mu},{\bold q}_{\nu}],
$$
or, more precisely, the closure of the commutator.  We want the $Q_{\mu\nu}$'s
to be self-adjoint operators which commute with each other and with the ${\bold
q}_{\mu}$'s.  The $Q_{\mu\nu}$'s, while unbounded, will together ``generate'' a
commutative $C^*$-algebra, their joint spectrum, and will correspond to
unbounded real-valued continuous functions on this spectrum.  Given a point,
$s$, of the spectrum, each $Q_{\mu\nu}$ can be evaluated at $s$, and the
resulting matrix, $Q_{\mu\nu}(s)$, will be a $4 \times 4$ matrix which is
skew-symmetric because of how the $Q_{\mu\nu}$'s are defined.  We require that
this matrix be in $\Sigma$.  In fact we require (as in \cite{4}) that $\Sigma$
be the joint spectrum of the $Q_{\mu\nu}$'s, so that each $Q_{\mu\nu}$ is
identified with the function on $\Sigma$ whose value at $\s \in \Sigma$ is just
$\s_{\mu\nu}$.

Finally, we require that all of this have an integrated Weyl form in terms of
one-parameter unitary groups generated by the ${\bold q}_{\mu}$'s (i\.e\. be a
``regular realization'' in the terminology of \cite{4}).  In fact, our
objective is to construct a $C^*$-algebra and specific affiliated elements such
that its non-degenerate representations give exactly all 
the regular realizations.

\heading
2.  The $C^*$-algebra
\endheading
To simplify notation we will work largely in a coordinate-free manner.  Let $V$
be any finite dimensional real vector space (eventually to be Minkowski space
$M_0$).  Let $V'$ denote its vector-space dual, so that each non-zero element
of $V'$ can be viewed as an unbounded function on $V$ (a coordinate function
for some basis of $V$).  Following the usage in \cite{4} we will sometimes
denote elements of $V$ by the letter $q$, and elements of $V'$ by the letters
$\a$ and $\b$, but we will also find it convenient to use $x$ or $v$ for
elements of $V$, and $p$ for elements of $V'$.

Let $GL(V)$ denote the group of invertible operators on $V$, and let $G$ denote
a closed subgroup of $GL(V)$ (eventually the Lorentz group for some Lorentz
metric on $V$).  Let $J$ be an invertible operator from $V'$ to $V$ which is
skew-symmetric in the sense that $J^t = -J$.  (Eventually $J$ will correspond
to the standard symplectic matrix $\s_0$ above.)  Now $GL(V)$ acts by
similarity on operators from $V'$ to $V$.  We let $\Sigma$ denote the orbit of
$J$ under $G$, that is, the set of all operators of the form $TJT^t$ for $T \in
G$.  Then $\Sigma$ is a homogeneous space for $G$, and is in particular a
manifold (perhaps not connected).

Our $C^*$-algebra will be a deformation of the commutative $C^*$-algebra 
\linebreak
$C_{\i}(\Sigma \times V)$ of continuous complex-valued functions on $\Sigma
\times V$ vanishing at infinity.  Specifically, for suitable functions in
$C_{\i}(\Sigma \times V)$ we want to define their deformed product by
$$
(f \times g)(\s,q) = \int_V \int_{V'} f(\s,q - \s(p))g(\s,q - v)e(v \cdot
p)dv\ dp,
$$
where $e$ denotes the function $e(t) = \exp(2\pi it)$, where $v \cdot p$
denotes the pairing between $V$ and $V'$ (which sometimes we will find more
convenient to write as $p(v)$), and where we have chosen a Haar (i\.e\.
Lebesgue) measure on $V$ and the corresponding Plancherel Lebesgue 
measure on $V'$.
This product is a general version of the Weyl, or Moyal, product.  We remark
that by including the $2\pi$ in our definition of $e$, we are, in effect,
choosing to follow the convention in \cite{5} rather than that in \cite{4}.

But the above definition of deformed product does not quite fit into 
the framework
of \cite{5} because, while the action of $V$ in \cite{5} is permitted to be
very complicated, $\s$ is supposed to be held fixed, while in the above
integral $\s$ is allowed to vary while the action is essentially just
translation.  One can certainly extend the general framework of \cite{5} to
cover this slightly different situation, but this note is not the place to do
that.  Instead we will show here that we can rearrange matters a little so that
the framework of \cite{5} does apply to yield the desired construction.  To do
this we will initially work on $G \times V$ rather than $\Sigma \times V$, and
then move to $\Sigma \times V$ near the end.

Actually, we will see that it is technically convenient for us to work more
generally on $E \times V$ where $E$ is any fixed open subset of $G$.  This is
because when we consider the unbounded operators affiliated with the
$C^*$-algebra we construct, it will be convenient for some steps to take as $E$
a bounded (open) subset of $G$.

For fixed $E$, define an action $\tau$ of $V$ on $E \times V$ by
$$
\tau_x(T,q) = (T,q+Tx).
$$
Let $\tau$ also denote the corresponding action on the $C^*$-algebra
$C_b(E\times V)$ of bounded continuous functions on $E \times V$.  This action
is not strongly continuous.  Let $B = C_u(E,V)$ denote the $C^*$-subalgebra of
all strongly continuous vectors for $\tau$, that is, functions $f$ such that $x
\mapsto \tau_xf$ is continuous on $V$ for the supremum norm on $C_b(E\times
V)$.  (The subscript ``$u$'' stands for ``uniformly continuous''.)  This puts
us directly in the framework of \cite{5,6}, so that we can construct
corresponding deformed $C^*$-algebras.  Let $B^{\i}$ denote the dense
$*$-subalgebra of smooth (i\.e\. infinitely differentiable) vectors in $B$.
Then for \linebreak 
$f$, $g \in B^{\i}$ we can define their deformed product (see page 83
of \cite{6}) by
$$
f\times_Jg = \int_{V'} \int_V \tau_{Jp}(f)\tau_v(g)e(v \cdot p)dv\ dp.
$$
The integrand will usually not be integrable, and so the integral must be
interpreted as an oscillatory integral by means of the factor $e(v\cdot p)$, as
discussed early in \cite{5}.  For any $(T,q) \in E\times V$ the above formula
gives, after using propositions $1.13$ or $2.12$ of \cite{5},
$$
(f\times_Jg)(T,q) = \int_{V'} \int_V f(T,q-TJT^tp)g(T,q-v)e(v\cdot p)dv\ dp.
$$
The natural involution is just complex conjugation, $f^*(T,q) =
(f(T,q))^-$.  We will denote by $B_J^{\i}$ the vector space $B^{\i}$ equipped
with the above deformed product and involution.  To view this explicitly as a
one-parameter deformation of the pointwise product, we should (as discussed in
\cite{5}) simply put a Planck's constant $\hslash$ in front of $J$, or more
appropriately in the present context, a Planck length $\l_p$.

Following \cite{5}, we put an operator norm on $B_J^{\i}$ as follows.  Let
${\Cal S}^B$ denote the vector-space of $B$-valued Schwartz functions on $V$,
with $B$-valued inner-product defined by
$$
\<\xi,\eta\>_B = \int_V (\xi(v))^*\eta(v)dv.
$$
Define a (real-valued) norm on ${\Cal S}^B$ by
$$
\|\xi\| = \|\<\xi,\xi\>_B\|^{1/2},
$$
where the norm on the right is that of $B$.  For $f \in B^{\i}$ define an
operator, $L_f$, on ${\Cal S}^B$ by
$$
(L_f\xi)(x) = \int_{V'} \int_V \tau_{x + Jp}(f)\xi(x+v)e(v\cdot p)dv\ dp.
$$
It is seen in theorem $4.1$ of \cite{Rfd} that $L_f$ is a bounded operator on
${\Cal S}^B$, that $f \mapsto L_f$ is a faithful $*$-representation of
$B_J^{\i}$ on
${\Cal S}^B$ for the $B$-valued inner-product, and that the corresponding norm
on $B_J^{\i}$ is a $C^*$-algebra norm.  We denote by $B_J$ the $C^*$-algebra
obtained by completing $B_J^{\i}$ for this norm.  It is represented on the
completion of ${\Cal S}^B$ for its norm.  Every state on $B$ can be composed
with the $B$-valued inner-product on ${\Cal S}^B$ to give an ordinary
inner-product on ${\Cal S}^B$, and hence a representation of $B_J$ on an
ordinary Hilbert space.  We obtain in this way a faithful family of ordinary
representations of $B_J$.

Let $A = C_{\i}(E \times V)$.  Then $A$ is an essential ideal in $B$ which is
carried into itself by $\tau$.  Then we have $A^{\i}$, $A_J^{\i}$, and $A_J$,
much as above.  Moreover, $A_J$ is an essential ideal in $B_J$ by proposition
$5.9$ of \cite{5}, so that we can view $B_J$ as a (unital) subalgebra of
the multiplier algebra, $M(A_J)$, of $A_J$.  It is $A_J$ which, for $V$
Minkowski space and $E = G$ the Lorentz group, will almost be our $C^*$-algebra
for the space-time uncertainty relations.  (We will still have to move the
situation to $\Sigma$.)

\heading
3.  The affiliated space-time operators
\endheading
Let $\a \in V'$.  We want to associate with $\a$ an unbounded operator, ${\bold
q}_{\a}$, affiliated with $A_J$ in the sense of Baaj \cite{1} and Woronowicz
\cite{8}.  The $C^*$-algebra ``generated'' by ${\bold q}_{\a}$ should be
isomorphic to $C_{\i}({\Bbb R})$, and should consist simply of the range of the
map from $C_{\i}({\Bbb R})$ to functions on $E \times V$ which sends $\varphi
\in C_{\i}({\Bbb R})$ to $\varphi \circ \a$ independent of the $E$-variable.
But there is no action of $V$ on ${\Bbb R}$ for which this map is equivariant
(for $\tau)$.  Since we need equivariance in order to invoke the functoriality
of the construction of \cite{5}, we enlarge the domain as follows.  Consider
$C_b(E \times {\Bbb R})$, and define an action, $\rho^{\a}$, of $V$ on $C_b(E
\times {\Bbb R})$ by
$$
(\rho_x^{\a}\psi)(T,r) = \psi(T,r - \a(Tx)).
$$
Define a homomorphism, $\Phi^{\a}$, from $C_b(E \times {\Bbb R})$ into $C_b(E
\times V)$ by
$$
(\Phi^{\a}(\psi))(T,q) = \psi(T,\a(q)).
$$
Then it is easily verified that $\Phi^{\a}$ is equivariant for $\rho^{\a}$ and
$\tau$.  Much as above, let $C_u^{\a}(E \times {\Bbb R})$ denote the algebra of
strongly continuous vectors for $\rho^{\a}$, so that its deformation
quantization, $(C_u^{\a}(E \times {\Bbb R}))_J$, is defined.  Then $\Phi^{\a}$
carries $C_u^{\a}(E \times {\Bbb R})$ into $C_u(E \times V)$, and so by the
functoriality of our deformation quantization construction (theorem $5.7$ of
\cite{5}) $\Phi^{\a}$ determines a $*$-homomorphism, $\Phi_J^{\a}$, from
$(C_u^{\a}(E \times {\Bbb R}))_J$ to $B_J$.

Let us now calculate the deformed product on $C_u^{\a}(E \times {\Bbb R})_J$,
and see that the product is in fact unchanged---no surprise in view of the
basically one-dimensional nature of the action $\rho^{\a}$.  If $\psi_1,\psi_2$
are $\rho^{\a}$-smooth vectors in \linebreak 
$C_u(E \times {\Bbb R})$, then much as for
$B^{\i}$ earlier
$$
(\psi_1 \times_J \psi_2)(T,r) = \int_{V'} \int_V \psi_1(T,r -
\a(TJT^tp))\psi_2(T,r - \a(v))e(v\cdot p)dv\ dp.
$$
But this ``sees'' only the $\a$-component of $v \in V$, and so by proposition
$1.3$ of \cite{6}, which is just a reinterpretation of proposition $1.11$ of
\cite{5}, we can reduce the domain of integration as follows.  Let $W$ denote
the kernel of $\a$ in $V$, so that $W^{\bot}$ is just the linear span of $\a$
in $V'$.  Then the above integral becomes
$$
\int_{W^{\bot}} \int_{V/W} \psi_1(T,r - \a(TJT^tp))\psi_2(T,r - \a(v))e(v\cdot
p)dv\ dp.
$$
But since $p$ is now running over the span of $\a$ and $TJT^t$ is
skew-symmetric, $\a(TJT^tp) = 0$.  From corollary $1.12$ of \cite{5} (basically
the Fourier inversion formula) it follows that
$$
\psi_1 \times_J \psi_2 = \psi_1\psi_2,
$$
the pointwise product.  The $C^*$-completion is then evidently $C_u(E \times
{\Bbb R})$ itself.

Thus we see that $\Phi^{\a}$ determines a homomorphism, $\Phi_J^{\a}$, of the
commutative $C^*$-algebra $C_u(E \times {\Bbb R})$ into $B_J$.  Since
$\Phi^{\a}$ for the undeformed algebras is injective, it follows from
proposition $5.8$ of \cite{5} that $\Phi_J^{\a}$ is injective.

Notice, however, that $\Phi^{\a}$ does not carry $C_{\i}(E \times {\Bbb R})$
into $C_{\i}(E \times V)$ (unless $V$ is one-dimensional) since functions in
its range are constant on the kernel of $\a$.  But $\Phi^{\a}$ will give a
``morphism'' from $C_{\i}(E \times {\Bbb R})$ to $C_{\i}(E \times V)$ in the
sense of Woronowicz \cite{8}.  This means that $\Phi^{\a}$ carries $C_{\i}(E
\times {\Bbb R})$ into the multiplier algebra, $C_b(E \times V)$, of $C_{\i}(E
\times V)$, and that $\Phi^{\a}(C_{\i}(E \times {\Bbb R}))C_{\i}(E \times V)$
(linear span) is dense in $C_{\i}(E \times V)$.  By the functoriality of the
deformation quantization construction with respect to morphisms, given in
theorem $3.1$ of \cite{7}, it follows that when we restrict $\Phi_J^{\a}$ to
$C_{\i}(E \times {\Bbb R})$ it determines a morphism to $A_J$.

We are now in a position to define the operators ${\bold q}_{\a}$ (which will
be the space-time operators when $V$ is Minkowski space).  Given $\varphi \in
C_{\i}({\Bbb R})$, view it as the corresponding function on $E \times {\Bbb R}$
independent of the $E$-variable.  Although the resulting function need not be
in $C_{\i}(E \times {\Bbb R})$, we clearly obtain in this way a morphism from
$C_{\i}({\Bbb R})$ to $C_{\i}(E \times {\Bbb R})$.  But morphisms can be
composed---see the comments near the end of section 0 of \cite{W}.  So if we
compose the above morphism with $\Phi_J^{\a}$, we obtain a morphism,
$\Theta_J^{\a}$, from $C_{\i}({\Bbb R})$ to $A_J$.  The unbounded function
$\iota(t) = t$ on ${\Bbb R}$ defines a self-adjoint operator affiliated with
$C_{\i}({\Bbb R})$, in the sense of \cite{W}.  But according to theorem $1.2$
of \cite{W} affiliated operators are transported by morphisms, and so there is
a well-defined self-adjoint operator, $\Theta_J^{\a}(\iota)$, affiliated with
$A_J$.  We set
$$
{\bold q}_{\a} = \Theta_J^{\a}(\iota).
$$

The one-parameter unitary group, $v_t$, generated by the operator $\iota$
affiliated to $C_{\i}({\Bbb R})$ is clearly $v_t(r) = e(tr)$ (up to conventions
about $2\pi$).  Under the morphism $\Theta_J^{\a}$ this is carried to the
one-parameter unitary group, $u_t^{\a}$, generated by ${\bold q}_{\a}$.  Since
we can scale $\a$, we only need $u_1^{\a} = e({\bold q}_{\a})$, which we denote
by $u_{\a}$.  We will show that the $u_{\a}$'s satisfy the appropriate Weyl
relations.

All the above is very natural, but there is an interesting technical point
which now needs to be emphasized, and which is responsible for our working in
the generality of the subsets $E$ of $G$.  Suppose $\varphi \in C_{\i}({\Bbb
R})$ is given, and view it as a function on $E \times {\Bbb R}$ as done above.
Then this function need not be in $C_u^{\a}(E \times {\Bbb R})$, because for $x
\in V$ we have
$$
\|\rho_x^{\a}(\varphi) - \varphi\| = \sup\{|\varphi(r-\a(Tx)) - \varphi(r)|: T
\in E \text{ and } r \in {\Bbb R}\},
$$
which need not go to $0$ as $x \to 0$, since the set $E$ over which $T$ varies
may be an unbounded subset of $GL(V)$.  Then even if $\varphi$ is a smooth
vector in $C_{\i}({\Bbb R})$ for translation, $\Theta^{\a}(\varphi)$ need not
be a strongly continuous vector in $C_b(E \times V)$, and so its deformed
product with elements of $A^{\i}$ can not be computed by the oscillatory
integral used earlier.  A very similar situation occurs in \cite{7}---see the
comments after definition $3.4$ of \cite{7}.  These include comments about how
to compute the deformed product by an approximation process, and analogous
comments apply here.

However, a simple argument shows that if $E$ is a {\it bounded} subset of
$GL(V)$ then indeed $\varphi$ viewed on $E \times {\Bbb R}$ is in $C_u(E \times
{\Bbb R})$, and moreover 
that if $\varphi$ is a smooth vector in $C_{\i}({\Bbb R})$ for
translation then it is a smooth vector also in $C_u(E \times {\Bbb R})$ for
$\rho^{\a}$, and thus $\Theta^{\a}(\varphi)$ will be in $B^{\i}$.  We need to
use this fact in order to check the appropriate Weyl relations for the
$u_{\a}$'s.  (A simple calculation at the level of functions, given below,
shows what these relations should be.  But we must continue to work in the
framework of \cite{5} in order to see that all is compatible with respect to
the operator norms which are present.)

Since we will now vary $E$, we need to make $E$ explicit in our notation.  So
we will now reserve $A$ for $C_{\i}(G \times V)$, and write ${}^EA$ for
$C_{\i}(E \times V)$, with  similar notation for $B$.  We have the evident
restriction map, $R^E$, from $C_{\i}(G \times V)$ to $C_b(E \times V)$, which
is a morphism from $A$ to ${}^EA$.  Since $E$ is open, the range of $R^E$ will
contain ${}^EA$, but it will usually contain much more.  Clearly $R^E$ is
equivariant for $\tau$ acting on both $A$ and ${}^EA$.  Thus $R^E$ determines a
morphism, $R_J^E$, from $A_J$ to ${}^EA_J$, by theorem $3.1$ of \cite{7}.

Since $E$ is open, each function in $C_{\i}(E \times V)$ can be viewed as a
function in $C_{\i}(G \times V)$ by setting it equal to $0$ off $E \times V$.
In this way ${}^EA$ is an ideal in $A$, carried into itself by $\tau$.  Thus
${}^EA_J$ is an ideal in $A_J$ by proposition $5.9$ of \cite{5}.  In
particular, for $f \in {}^EA^{\i}$, viewed as in $A^{\i}$, we can consider
$R_J^E(u_{\a} \times_J f) = R_J^E(u_{\a}) \times_J R_J^E(f)$.  Now if $E$ is
bounded, $R^E(u_{\a})$ is a smooth vector in ${}^EB$.  This is exactly our
reason for restricting to bounded $E$.  For then we can calculate the deformed
product directly by our oscillatory integrals.  To lighten our notation we will
omit the $R^E$ but work on $E$.  Then we have:

\proclaim{3.1 Lemma}
For bounded $E$, for $f \in {}^EA^{\i}$, and for $\a \in V'$, we have
$$
(u_{\a} \times_J f)(T,q) = e(q \cdot \a)f(T,q + TJT^t\a)
$$
for $(T,q) \in G \times V$, in the sense that $u_{\a} \times_J f$ as element of
${}^EA_J$ is represented by the function on the right.
\endproclaim

\demo{Proof}
For $(T,q) \in E \times V$ we have
$$
\align
(u_{\a} \times_J f)(T,q) &= \int_{V'} \int_V e(\a(q-TJT^tp)f(T,q-v)e(v\cdot
p)dv\ dp \\
&= e(q \cdot \a) \int_{V'} \int_V e((v+TJT^t\a) \cdot p)f(T,q - v)dv\ dp \\
&= e(q \cdot \a)f(T,q + TJT^t\a),
\endalign
$$
where we have used the Fourier inversion formula in the last step.  Thus
$R_J^E(u_{\a} \times_J f)$ is given by the function on the right.  But $R_J^E$
is injective by proposition $5.8$ of \cite{5}.  Q\.E\.D\.
\enddemo

\proclaim{3.2 Theorem}
For $\a,\b \in V'$ we have
$$
u_{\a} \times_J u_{\b} = e(Q_{\a\b})u_{\a+\b}
$$
in $M(A_J)$, where $Q_{\a\b}$ is defined on $G \times V$ by
$$
Q_{\a\b}(T,q) = \b(TJT^t\a).
$$
\endproclaim

\demo{Proof}
We remark first that $e(Q_{\a\b})$ is invariant under $\tau$ since $Q_{\a\b}$
is independent of $q$.  Thus $e(Q_{\a\b})$ is a smooth vector in $B$, whose
deformed product with any function in $B^{\i}$ is just the pointwise product
(by corollary $2.13$ of \cite{5}).  (In particular, all the $e(Q_{\a\b})$'s
commute among themselves for the deformed product, as do thus the
$Q_{\a\b}$'s.) In this way $e(Q_{\a\b})u_{\a+\b}$ has
meaning.  Now for any bounded $E$ and for any $f \in {}^EA^{\i}$ we have
(omitting $R^E$ from our notation)
$$
\align
(u_{\a} \times_J (u_{\b} \times_J f))&(T,q) = e(q \cdot \a)(u_{\b} \times_J
f)(T,q + TJT^t\a) \\
&= e(q\cdot \a)e((q + TJT^t\a) \cdot \b)f(T,q + TJT^t\a + TJT^t\b) \\
&= e(\b(TJT^t\a))(u_{\a+\b} \times_J f)(T,q),
\endalign
$$
for $(T,q) \in E \times V$.  But ${}^EA_J$ is an ideal in $A_J$, and so, using
associativity, we have
$$
(u_{\a} \times_J u_{\b}) \times_J f = (e(Q_{\a\b})u_{\a+\b}) \times_J f.
$$
But the union of the ${}^EA$'s is a dense $\tau$-invariant ideal in $A$.  It is
easily seen from the proof of proposition $2.17$ of \cite{5} that the smooth
elements of this dense ideal are dense in $A_J$.  Consequently $u_{\a} \times_J
u_{\b} = e(Q_{\a\b})u_{\a+\b}$ as desired.  \qed
\enddemo

This theorem shows that the unbounded operators ${\bold q}_{\a}$ affiliated
with $A_J$ have an integrated Weyl form, so that they give a ``regular''
realization in the sense of \cite{4}, satisfying
$$
[{\bold q}_{\a},  {\bold q}_{\b}] = Q_{\a\b}.
$$
The operators $Q_{\a\b}$ are affiliated in the evident way with $C_{\i}(G)$,
and so with $A_J$ through the evident morphism of $C_{\i}(G)$ to $C_{\i}(G
\times V)$.  For the case considered in \cite{4} in which $V$ is Minkowski
space, $G$ is the Lorentz group, and $J$ is the standard symplectic matrix as
in $3.28$ of \cite{4}, the above operators will satisfy the conditions desired
in \cite{4}.

\heading
4.  The space-time algebra
\endheading
In the situation discussed in the previous section we notice that the
dependence of $Q_{\a\b}$ on $T$ is really only on $TJT^t$, and that in Lemma~2
the $T$-dependence of the product with $u_{\a}$ is only on $TJT^t$.  But the
$TJT^t$'s comprise the orbit, $\Sigma$, of $J$ under $G$.  This indicates that
our algebra $C_{\i}(G \times V)$ is unnecessarily large, and that we should be
working on $\Sigma \times V$.  We did not do this initially because there is no
convenient action of $V$ 
there which we could use so as to apply the constructions of
\cite{5}.  But there is a simple device by which we can now pass to $\Sigma
\times V$.

Let $H$ denote the stability subgroup of $J$ in $G$, consisting of those $S \in
G$ such that $SJS^t = J$  (so that $\Sigma$ is identified with $G/H$).  
Since the set of $S$'s in $\text{Aut}(V)$ which
satisfy this relation is just the symplectic group for $J$, we see that $H$ is
just the intersection of $G$ with this symplectic group.  Define a (free and
proper) right action, $\g$, of $H$ on $G \times V$ by
$$
\g_S(T,q) = (TS,q),
$$
and let $\g$ also denote the corresponding action on $B = C_u(G \times V)$.
This action does not commute with $\tau$ in general.  But for $f \in B$ we have
$$
\align
(\tau_x(\g_Sf))(T,q) &= f(TS^{-1},q - Tx) \\
&= (\g_S(\tau_{Sx}f))(T,q),
\endalign
$$
that is, $\tau_x\g_S = \g_S\tau_{Sx}$.  But this is exactly the condition
(together with $(SJS^t = J)$ which according to proposition $10.4$ of \cite{5}
insures that $\g_S$ determines an automorphism of $B_J$.  Thus $\g$ gives a
right action (which we still denote by $\g$) of $H$ on $B_J$.  
We do not need to
explore here the strong continuity of this action, since what we are interested
in is the fixed-point subalgebra, $D_J$.  Thus for this purpose we can treat
$H$ as a discrete group.

At the function level it is clear that there are many elements in $D_J$.  Any
$f \in C_c(\Sigma \times V)$ which is smooth in the $V$-direction will lift
to
an element of $D_J$.  It is clear from Lemma $3.1$ that taking products with
any $u_{\a}$ will carry $D_J$ into itself, so that the $u_{\a}$'s can be viewed
as multipliers of $D_J$.  They will still satisfy the Weyl relations
$$
u_{\a} \times_J u_{\b} = e(Q_{\a\b})u_{\a+\b},
$$
but now we can view $Q_{\a\b}$ as the function on $\Sigma$ (or $\Sigma \times
V$) defined by
$$
Q_{\a\b}(\s) = \b(\s(\a)).
$$
This means that we still have the unbounded operators ${\bold q}_{\a}$, but now
affiliated with $D_J$.  They still satisfy
$$
[{\bold q}_{\a},{\bold q}_{\b}] = Q_{\a\b}
$$
as above. It is easily
seen that we obtain:

\proclaim{4.1 Theorem}
The $C^*$-algebra $D_J$ together with its 
affiliated unbounded operators ${\bold q}_{\a}$ and 
corresponding unitary operators $u_{\a}$, is the universal $C^*$-algebra
for ``regular realizations'' of the commutation relations over $\Sigma$, that
is, integrable representations for which
$$
[{\bold q}_{\a},{\bold q}_{\b}] = Q_{\a\b},
$$
where the $Q_{\a\b}$'s commute among themselves 
and with the ${\bold q}_{\a}$'s, and
for which the joint spectrum of the $Q_{\a\b}$'s is a subset of $\Sigma$.
\endproclaim

When $V$ is Minkowski space, etc, then $D_J$ is essentially the
$C^*$-algebra of \cite{4} for the space-time uncertainty relations.

We could continue by applying other parts of \cite{5}, for example chapter~8 to
discuss the continuous field aspect which is discussed in \cite{4}, or
chapter~9 to discuss the semi-classical limit.  But this is all relatively
straightforward, and so we will stop our discussion here.

\heading
References
\endheading

\item{[1]} S.~Baaj:  {\it Multiplicateurs non born\'es}, Th\`ese 3eme cycle,
Universit\'e Paris VI, 1980.

\item{[2]} S.~Doplicher:  {\it Quantum physics, classical gravity and
non-commutative spacetime},  Proc. XIth International Congress of
Mathematical Physics, Paris, 1994, International Press,
Cambridge Ma, 324-329.

\item{[3]} S.~Doplicher, K.~Fredenhagen, and J.~E.~Roberts:  {\it Spacetime
quantization induced by classical gravity}, Phys. Lett. B, {\bf 331} (1994),
39-44.

\item{[4]} S.~Doplicher, K.~Fredenhagen, and J.~E.~Roberts:  {\it The quantum
structure of spacetime at the Planck scale and quantum fields}, Commun. Math.
Phys., 172 (1995), 187-220.

\item{[5]} M.~A.~Rieffel:  {\it Deformation Quantization for actions of ${\Bbb
R}^d$}, Memoirs A\.M\.S\. {\bf 506}, Amer. Math. Soc., Providence, 1993.

\item{[6]} M.~A.~Rieffel:  {\it Quantization and $C^*$-algebras}, Contemp.
Math., {\bf 167} (1994), 67-97.

\item{[7]} M.~A.~Rieffel:  {\it Non-compact quantum groups associated with
Abelian subgroups}, Commun. Math. Phys., {\bf 171} (1995), 181-201.

\item{[8]} S.~L.~Woronowicz:  {\it Unbounded elements affiliated with
$C^*$-algebras and non-compact quantum groups}, Commun. Math. Phys., {\bf 136}
(1991), 399-432.

\enddocument
\bye